\documentclass[reqno]{amsart}

\input xy
\xyoption{all}
\usepackage{epsfig}
\usepackage{color}
\usepackage{amsthm}
\usepackage{amssymb}
\usepackage{amsmath}
\usepackage{amscd}
\usepackage{amsopn}
\usepackage{url}
\usepackage{hyperref}\hypersetup{colorlinks}

\usepackage{mathabx}

\usepackage{eufrak}

\usepackage{color} 

\definecolor{darkred}{rgb}{1,0,0} 
\definecolor{darkgreen}{rgb}{0,0.8,0}
\definecolor{darkblue}{rgb}{0,0,1}

\hypersetup{colorlinks,
linkcolor=darkblue,
filecolor=darkgreen,
urlcolor=darkred,
citecolor=darkgreen}

\numberwithin{equation}{section}

\newcommand{\labell}[1] {\label{#1}}

\numberwithin{equation}{section}
\newtheorem {Theorem}{Theorem}
\numberwithin{Theorem}{section}

\newtheorem {Lemma}[Theorem]    {Lemma}

\theoremstyle{definition}

\theoremstyle{remark}
\newtheorem{Remark}[Theorem]{Remark}
\newtheorem{Example}[Theorem]{Example}

\expandafter\chardef\csname pre amssym.def at\endcsname=\the\catcode`\@
\catcode`\@=11
\def\undefine#1{\let#1\undefined}
\def\newsymbol#1#2#3#4#5{\let\next@\relax
 \ifnum#2=\@ne\let\next@\msafam@\else
 \ifnum#2=\tw@\let\next@\msbfam@\fi\fi
 \mathchardef#1="#3\next@#4#5}
\def\mathhexbox@#1#2#3{\relax
 \ifmmode\mathpalette{}{\m@th\mathchar"#1#2#3}%
 \else\leavevmode\hbox{$\m@th\mathchar"#1#2#3$}\fi}
\def\hexnumber@#1{\ifcase#1 0\or 1\or 2\or 3\or 4\or 5\or 6\or 7\or 8\or
 9\or A\or B\or C\or D\or E\or F\fi}

\font\teneufm=eufm10
\font\seveneufm=eufm7
\font\fiveeufm=eufm5
\newfam\eufmfam
\textfont\eufmfam=\teneufm
\scriptfont\eufmfam=\seveneufm
\scriptscriptfont\eufmfam=\fiveeufm

\catcode`\@=\csname pre amssym.def at\endcsname

\newcommand{\EE}{{\mathcal E}}
\newcommand{\FF}{{\mathcal F}}

\newcommand{\hf}{\hat{f}}

\newcommand{\A}{{\mathcal A}}

\newcommand{\PP}{{\mathcal P}}

\def    \F      {{\mathbb F}}

\def    \C      {{\mathbb C}}
\def    \R      {{\mathbb R}}

\def    \Z      {{\mathbb Z}}

\def    \T      {{\mathbb T}}

\def    \12    {{\frac{1}{2}}}

\def    \SB     {\operatorname{SB}}

\def    \SU     {\operatorname{SU}}

\def    \GL    {\operatorname{GL}}

\def    \H     {\operatorname{H}}

\def    \CL     {\operatorname{CL}}

\begin{document}


\setlength{\smallskipamount}{6pt}
\setlength{\medskipamount}{10pt}
\setlength{\bigskipamount}{16pt}





\title[Hyperk\"ahler Arnold Conjecture and its Generalizations]{Hyperk\"ahler 
Arnold Conjecture and its Generalizations}

\author[Viktor Ginzburg]{Viktor L. Ginzburg}
\author[Doris Hein]{Doris Hein}

\address{Department of Mathematics, UC Santa Cruz,
Santa Cruz, CA 95064, USA}
\email{ginzburg@math.ucsc.edu}
\email{dhein@ucsc.edu}

\date{\today} \thanks{The work is partially supported by the NSF and
by the faculty research funds of the University of California, Santa
Cruz.}


\begin{abstract} We generalize and refine the hyperk\"ahler Arnold
  conjecture, which was originally established, in the non-degenerate
  case, for three-dimensional time by Hohloch, Noetzel and Salamon by
  means of hyperk\"ahler Floer theory. In particular, we prove the
  conjecture in the case where the time manifold is a multidimensional
  torus and also establish the degenerate version of the
  conjecture. Our method relies on Morse theory for generating
  functions and a finite-dimensional reduction along the lines of the
  Conley--Zehnder proof of the Arnold conjecture for the torus.
\end{abstract}

\maketitle

\tableofcontents

\section{Introduction}
\labell{sec:intro}
The main objective of this paper is to prove a generalization of the
hyperk\"ahler Arnold conjecture originally established via
hyperk\"ahler Floer theory by Hohloch, Noetzel and Salamon in
\cite{HNS}.

The setting of the hyperk\"ahler Arnold conjecture is similar to its
standard Hamiltonian counterpart, but the time manifold is
three-dimensional ($\T^3$ or $\SU(2)$ rather than $S^1$) and the
target manifold is equipped with a hyperk\"ahler rather than a
symplectic structure. The space of maps from the time manifold to the
target manifold carries a suitably defined action functional, akin to
the standard action functional in Hamiltonian mechanics, provided that
a version of a Hamiltonian is also furnished. In the spirit of the
Arnold conjecture, the main result of hyperk\"ahler Floer theory
developed \cite{HNS} is that the number of critical points of the
action functional is bounded from below by the sum of Betti numbers of
the target manifold whenever the action functional is Morse.  For
technical reasons, the target manifold must be flat.

Our main goal is to show that this version of the Arnold conjecture
can be further generalized and refined.  We prove an analog of the
conjecture (both the degenerate and non-degenerate case) for the time
manifold $\T^r$ and a target space equipped with $r$ flat
``anti-commuting'' K\"ahler structures. More precisely, the target
space is a compact quotient of a representation of a Clifford algebra.
In the degenerate case, the lower bound is given in terms of the
cup-length of the target space. We also prove a version of the
degenerate Arnold conjecture for the time manifold $\SU(2)$ and a
flat hyperk\"ahler target space.

In contrast with \cite{HNS}, the argument we utilize to prove these
results is not precisely Floer theoretic, but rather it is a
finite-dimensional approximation combined with Morse or
Ljusternik--Schnirelman theory for generating functions, following the
line of reasoning from \cite{CZ}. The difference is, from our
perspective, rather technical and the two methods usually give 
the same results when they both apply, with, perhaps, the
finite-dimensional approximation approach having a slight edge. (Of
course, in the context of Hamiltonian dynamics, Floer theory has a
much broader range.)

\subsection{Acknowledgments} We are grateful to Sonja Hohloch, Richard
Montgomery, Gregor Noetzel, Jie Qing, Dietmar Salamon, Alan Weinstein,
and Martin Weissman for useful discussions.

\section{Main Results}
\label{sec:main}
Let $V$ be a vector space equipped with $r$ symplectic structures
$\omega_1,\ldots,\omega_r$, which are all compatible with the same
inner product $\left<\,,\right>$. In other words, there exist
orthogonal (with respect to $\left<\,,\right>$) operators $J_1,\ldots,
J_r$ on $V$ such that $J_l^2=-I$ for all $l$, i.e., these operators
are complex structures, and
$$
\left<X,Y\right>=\omega_l(X,J_l Y) \text{ for all $X$ and $Y$ in $V$}.
$$
Assume furthermore that the complex structures $J_l$ anti-commute:
\begin{equation}
\label{eq:anti-comm}
J_lJ_j+J_jJ_l=0 \text{ whenever $l\neq j$}.
\end{equation}

Such a collection of complex (or equivalently symplectic) structures
can exist for arbitrarily large values of $r$, depending on the dimension of $V$. 
It exists if and only if the unit sphere in $V$ admits $r$
linearly independent vector fields; see \cite[Chapter 12 and 16]{Hu}
and, in particular, pp.\ 152--154 therein. More specifically, let
$\dim V= 2^{4d+c}b$, where $d\geq 0$ and $0\leq c\leq 3$ are integers
and $b$ is odd. Then the maximal value of $r$ for $V$ is
$8d+2^c-1$. In fact, equipping $V$ with the structures
$J_1,\ldots,J_r$ is equivalent to turning $V$ into an
(orthogonal) representation of the Clifford algebra of a negative
definite quadratic form on $\R^r$. Note also that the forms $\omega_l$
generate a ``pencil'' of symplectic structures, i.e., as is easy to
see, any non-trivial linear combination $\omega=\sum
\lambda_l\omega_l$ is symplectic. Likewise, a linear combination
$J=\sum \lambda_l J_l$ is, up to a factor, a complex structure. More
precisely, $J^2=-(\sum \lambda_l^2) I$.

\begin{Example}[Hyperk\"ahler structures]
A standard example of a vector space with such structures is a
hyperk\"ahler vector space. In this case, $r=3$ and the complex
structures $J_l$ satisfy the quaternionic relations, i.e., in addition
to \eqref{eq:anti-comm} we also have $J_1J_2=J_3$.
\end{Example}

Let now $W$ be a smooth compact quotient of $V$ by a group of
transformations preserving all of the above structures on $W$. For
instance, $W$ can be the quotient of $V$ by a lattice. (There are,
however, other examples; see, e.g., \cite[p.\ 2548]{HNS}.)

Furthermore, let us fix a closed manifold $M$ equipped with a volume
form $\mu$ and $r$ divergence--free vector fields
$v_1,\ldots,v_r$. This manifold will take the role of ``time'' in
Hamiltonian dynamics.  More specifically, the following two examples
are of interest to us.

\begin{Example}[The torus]
\label{ex:torus}
In this example, $M$ is the $r$-dimensional torus $\T^r=\R^r/\Z^r$
with angular coordinates $t_1,\ldots, t_r$, the vector fields $v_l$
are the coordinate vector fields $\partial_{t_l}$, and
$\mu=dt_1\wedge\ldots\wedge dt_r$. More generally, we can replace the
coordinate vector fields by any basis of vector fields with constant
coefficients.
\end{Example}

\begin{Example}[The special unitary group $\SU(2)$]
\label{ex:su(2)}
Here $r=3$ and $M=\SU(2)$ is equipped with the (probability) Haar measure
$\mu$. The vector fields $v_l$ are the right-invariant
vector fields whose values at the unit $e$ are:
\begin{equation}
\label{vectorfields}
v_1(e) =
\left(
\begin{array}{rr}
0 & i \\
i & 0 \\
\end{array}
\right)
,
\quad
v_2(e) =
\left(
\begin{array}{rr}
0 &-1 \\
1 & 0 \\
\end{array}
\right)
,
\quad
v_3(e) =
\left(
\begin{array}{rr}
i & 0 \\
0 & -i \\
\end{array}
\right).
\end{equation}
More generally, we may replace $\SU(2)$ by the homogeneous space
$M=\SU(2)/G$, where $G \subset \SU(2)$ is a discrete subgroup.
The vector fields $v_l$ naturally descend to this quotient.
\end{Example}

By analogy with Hamiltonian dynamics, a Hamiltonian is a smooth function
$$
H\colon M\times W\to \R.
$$
The action functional $\A_H$ is defined on the space $\EE$ of
$C^\infty$-smooth (or just $C^2$), null-homotopic maps $f\colon M\to
W$. We introduce $\A_H$ in two steps. First, let $F\colon [0,\, 1]\times M\to
W$ be a homotopy between $f$ and the constant map. This is an analog
of a capping in the definition of the standard Hamiltonian action
functional. The unperturbed action functional is 
$$
\A(f)=-\sum_l \int_{[0,\, 1]\times M} F^* \omega_l \wedge i_{v_l}\mu.
$$
It is routine to check that $\A(f)$ is well-defined, i.e., independent
of $F$. (Here it would be sufficient to assume that, e.g., the universal covering of $W$
is contractible.) Finally, the total or perturbed action functional is
\begin{equation}
\label{eq:action}
\A_H(f)=\A(f)-\int_M H(f)\mu.
\end{equation}
For instance, when $r=1$ and $M=\T^1$, we obtain the ordinary action
functional of Hamiltonian dynamics. Furthermore, it is easy to see
that in the setting of Example \ref{ex:torus} (with $r=3$) or of Example
\ref{ex:su(2)} the perturbed and unperturbed action functionals
coincide, up to a sign, with those defined in \cite{HNS}.

The differential of $\A$ at $f\in\EE$ is 
$$
(d\A)_f (w)=\sum_l\int_M \omega_l(L_{v_l}f,w)\mu,
$$
where $w\in T_f\EE$ is a vector field along $f$. Thus, the gradient of
$\A$ with respect to the natural $L^2$-metric on $\EE$ is a Dirac type operator
$$
\nabla_{L^2} \A (f) =\sum_l J_l L_{v_l} f =: \partialslash f.
$$
Hence, we have
$$
\nabla_{L^2} \A_H (f) = \partialslash f -\nabla H (f),
$$
where $\nabla H$ denotes the gradient of $H$ along $W$.  As a
consequence, the critical points of $\A_H$ are solutions $f\in\EE$ of
the equation
\begin{equation}
\label{eq:crit-pnts}
\partialslash f =\nabla H (f).
\end{equation}

At a critical point $f$ of $\A_H$, the Hessian $d^2_f \A_H$ is defined in
the standard way as the second variation of $\A_H$. This is a quadratic
form on $T_f\EE$ equal to the $L^2$-pairing with the linearization of
$\nabla_{L^2} \A_H$ at $f$.  We call $f$ a non-degenerate critical
point when this operator $T_f\EE\to T_f\EE$ is one-to-one, cf.\
\cite[p.\ 2559]{HNS}. A Hamiltonian $H$ is said to be non-degenerate
when all critical points of $\A_H$ are non-degenerate. In the setting of
Examples \ref{ex:torus} and \ref{ex:su(2)}, non-degeneracy is a
generic condition on $H$, i.e., the set of non-degenerate Hamiltonians
is residual in $C^\infty(M\times W)$. (The proof in \cite[p.\
2574--2576]{HNS} covers Example \ref{ex:su(2)} and carries
over to Example \ref{ex:torus} for all $r$ with straightforward
modifications.)

Finally, denote by $\CL(W)$ the cup-length of $W$, i.e., the maximal
number of elements in $\H_{*>0}(W;\F)$ such that their cup-product is
not equal to zero, also maximized over all fields $\F$.  Likewise, let
$\SB(W)$ (the sum of Betti numbers) stand for $\sum_j \dim_\F \H_j(W;\F)$,
maximized again over all $\F$.

In the spirit of the Arnold conjecture and of \cite{HNS}, our main
result is

\begin{Theorem}
\label{thm:main}
Assume that $M$ is as in Example \ref{ex:torus}, or that $V$ is
hyperk\"ahler and $M$ is as in Example \ref{ex:su(2)}. Then for any
Hamiltonian $H$, the action functional $\A_H$ has at least $\CL(W)+1$
critical points. If $H$ is non-degenerate, the number of critical
points is bounded from below by $\SB(W)$.
\end{Theorem}

We emphasize that the non-degenerate case of this theorem  was
originally proved in \cite{HNS} in the setting of a hyperk\"ahler target 
space and the domain being either $M=\SU(2)$ or $M=\T^3$.

Theorem \ref{thm:main} suggests that in this context a version of
Hamiltonian Floer theory can be developed beyond the setting where the
target space $W$ is hyperk\"ahler and the domain $M$ is hypercontact
as in \cite{HNS}. It appears that more generally a collection, as
above, of $r$ symplectic and complex structures on $W$ may be
sufficient for such a theory. Note however that manifolds equipped
with such structures must be extremely rare, cf. \cite[Chapter
21]{GHJ}. For instance, once $r\geq 2$, every such a manifold is
automatically hyperk\"ahler with the third complex structure
$J_1J_2$. The authors are not aware of any non-flat example where
$r>3$. Note also that similar, although not quite identical, types of
structures (at least on the complex side of the story) are considered
in \cite{MS,Jo}.  Pencils of symplectic structures also arise on the
point-wise (i.e., linear algebra) level on the manifolds equipped with
fat fiber bundles introduced in \cite{We:fat} or fat distributions;
see \cite[Section 5.6]{Mo} and references therein, and also
\cite{FZ}. It is less clear what in this setting the right structure
on the time manifold $M$ should be. We examine further generalizations
of the hyperk\"ahler Arnold conjecture elsewhere.

\begin{Remark}
\label{rmk:ellipticity}
  In the context of Floer theory, two properties of the operator
  $\partialslash$, hidden in our proof, are particularly
  important. Namely, the operator $\partialslash$ and the operator
  $\partial_s-\partialslash$ on $\R\times M$ must both be elliptic on
  the space of $V$-valued functions on $M$. To see when this is the
  case, let us assume for the sake of simplicity that the vector
  fields $v_l$ form a basis at every point of $M$. Then
  $\partialslash$ is elliptic if and only the symbol
  $\sigma(\partialslash)=\sum \lambda_l J_l$ is invertible for all
  non-zero (co)vectors $\lambda=(\lambda_1,\ldots, \lambda_r)$. This
  is clearly the case when, as above, the linear operators $J_l$ are
  anti-commuting complex structures; for then
  $\sigma(\partialslash)^2=-(\sum \lambda_l^2)I$. In a similar vein,
  $\partial_s-\partialslash$ is elliptic if and only if
  $\sigma(\partial_s-\partialslash)=\lambda_0I-\sum \lambda_l J_l$ is
  invertible for all $(\lambda_0,\lambda)\neq 0$. This is again
  automatically the case in our setting.
\end{Remark}

\begin{Remark}
\label{rmk:non-compact}
Theorem \ref{thm:main} extends to the case where the manifold $W$ is a
non-compact quotient of $V$ without any significant changes in the
proof. However, now certain restrictions must be imposed on the
behavior of the Hamiltonian $H$ at infinity and the lower bounds on
the number of critical points may possibly depend on these restrictions. To be more
specific, let us assume that a finite covering $W'$ of $W$ is a
Riemannian product of a flat torus and a Euclidean space $V'$. (For
instance, $W$ can be an iterated cotangent bundle of a flat 
manifold; it is not hard to see that this $W$ carries the required
structure.) Then it suffices to require the lift of $H$ to $M\times
W'$ to coincide outside a compact set with a non-degenerate quadratic
form on $V'$ with constant coefficients. In this case, the lower bounds on the
number of critical points are again  $\CL(W)+1$ and, respectively, $\SB(W)$.
\end{Remark}

\section{Proof of Theorem \ref{thm:main}}
\label{sec:proofs}

As has been pointed out in the introduction, the argument follows
closely the finite-dimensional reduction method of Conley and Zehnder, 
\cite{CZ}. The method utilizes the Fourier expansion of $f\colon M\to
W$ over $M$ to reduce the problem to now standard finite-dimensional
Morse theory for generating functions. In fact, when $M=\T^r$, the
proof carries over essentially word-for-word with hardly more than
notational changes. The case of $M=\SU(2)$ is more involved. For then
we use Fourier analysis on $\SU(2)$ -- the Peter--Weyl theorem
-- entailing somewhat lengthier calculations. In both cases, the main
point of the proof is obtaining an explicit expression for
$\partialslash f$ in terms of the Fourier expansion of $f$. Once this
is done, we faithfully adhere to the line of reasoning from \cite{CZ},
and hence omit here some straightforward, technical details of the
proof.

\subsection{The $\T^r$-case} 
\label{sec:torus}
Throughout the proof, we will assume that $v_1,\ldots,v_r$ are the
coordinate vector fields on $M=\T^r=\R^r/\Z^r$. The case of an
arbitrary basis of constant vector fields can be handled in a similar
way.

Furthermore, let us first assume that $W$ is the quotient of a vector
space $V$ by a lattice. (As a consequence, $W$ is a torus.) We will 
discuss the modifications needed to deal with the general case at the
end of the proof.
 
In what follows, it will be convenient to
view $V$ as a complex vector space, equipped with one of the complex
structures $J_l$, say, $J=J_r$. Since $f$ is null-homotopic, it can be
lifted to a map $\tilde{f}\colon M\to V$. Consider the Fourier
expansion of $\tilde{f}$:
$$
\tilde{f}(t)=\sum_k \exp(2\pi k\cdot t J) \hat{f}_k,
$$
where $t=(t_1,\ldots,t_r)\in \T^r$ and $k=(k_1,\ldots,k_r)\in\Z^r$ and
the Fourier coefficients $\hat{f}_k$ are elements of $V$. Note that
among these coefficients, the coefficients with $k\neq 0$ are
completely determined by $f$ and independent of the lift. (This is the
point where it is essential that $W$ is the quotient of $V$ by a
lattice.) The mean value $\hat{f}_0$ depends on the lift $\tilde{f}$,
but its image in $W$ is again completely determined by~$f$. Hence we
can, keeping the same notation $\hat{f}_0$ for the mean value,
unambiguously express $f$ as
\begin{equation}
\label{eq:fourier-exp}
f(t)=\sum_k \exp(2\pi k\cdot t J) \hat{f}_k,
\end{equation}
where now $\hat{f}_0\in W$ and $\hat{f}_k\in V$ when $k\neq 0$.

In other words, here we view $\EE$ as an infinite-dimensional vector
bundle over $W$ with projection map $f\mapsto \hat{f}_0$. This vector
bundle is trivial and its fiber $\FF$ is canonically isomorphic to the
space of smooth maps $M\to V$ with zero mean. The Fourier expansion
allows us, using self-explanatory notation, to regard $\EE$ as a
sub-bundle in $W\times L^2_0(M,V)$.

Our next goal is to obtain an explicit expression for $\partialslash
f$ in terms of the Fourier expansion \eqref{eq:fourier-exp}. As we
will soon see, the operator $\partialslash$ block-diagonalizes once we
group together the $k$th and $(-k)$th terms in
\eqref{eq:fourier-exp}. (Note that since $\partialslash$ kills
constant terms we can view it as either a linear operator on $\FF$ or
a fiberwise linear operator on $\EE=W\times \FF$ independent of the
point of the base.) To be more precise, let $k^*$ stand for a pair
$(-k,k)$, with $k\neq 0$. The pair is ordered lexicographically, i.e.,
so that the first non-zero component of $k$ is positive. 
Let $F_{k^*}$ be the subspace of $\FF$ formed by
functions $\exp(-2\pi k\cdot t J)X+\exp(2\pi k\cdot t J)Y$ with $X$
and $Y$ in $V$. Note that $L^2_0(M,V)$ is the $L^2$-direct
sum of the spaces $F_{k^*}$ for all pairs $k^*$. Below, we will use
the identification $F_{k^*}=V\oplus V$, where the first term
corresponds to $-k$ and the second one to $k$, and denote by $I$ the
identity operator on $V$.

\begin{Lemma}
\label{lemma:torus}
The space $F_{k^*}$ is invariant under $\partialslash$ and on this space,
$\partialslash$ acts as
\begin{equation}
\label{eq:A_k}
A_{k^*}=
2\pi 
\left[
\begin{array}{cc}
k_r I & -J\sum_{l=1}^{r-1} k_l J_l\\
J\sum_{l=1}^{r-1} k_l J_l & -k_r I
\end{array}
\right].
\end{equation}
Furthermore, $A_{k^*}$ is invertible and
\begin{equation}
\label{eq:A_k-inverse}
A_{k^*}^{-1}=
\frac{1}{4\pi^2\|k\|^2} A_{k^*} 
,
\end{equation}
where $\|k\|^2=k_1^2+\ldots+k_r^2$, and
\begin{equation}
\label{eq:A_k-inverse-norm}
\|A_{k^*}^{-1}\|=\frac{1}{2\pi\|k\|}.
\end{equation}
\end{Lemma}

\begin{Remark}
\label{rmk:lemma-torus} 
This lemma is more precise than is really necessary for the proof. In
fact, explicit expressions for $A_{k^*}$, its inverse and the norm
of the inverse are irrelevant. It would be sufficient to just know that
$A_{k^*}$ is invertible and that $\|A_{k^*}^{-1}\|=O(1/\|k\|)$.
\end{Remark}

\begin{proof}[Proof of the lemma]
Recall that $M$ is a torus, $v_l=\partial/\partial
t_l$ and
$$
\partialslash f= \sum J_l\frac{\partial f}{\partial t_l}.
$$
Thus, as a straightforward calculation shows,
$$
\partialslash f =2\pi\sum_ k\exp(2\pi k\cdot  t J)\left( J\sum_{l=1}^{r-1} 
k_l J_l \hf_{-k} - k_r \hf_k\right).
$$
Here we use the fact that $J=J_r$ anti-commutes with $J_l$ for
$l=1,\ldots, r-1$. This expression shows that $F_{k^*}$ is invariant
under $\partialslash$ and immediately implies \eqref{eq:A_k}.  Now
\eqref{eq:A_k-inverse} is straightforward to check using again the fact
that the complex structures $J_l$ with $l=1,\ldots,r$ anti-commute. To
finish the proof of the lemma, it remains to establish
\eqref{eq:A_k-inverse-norm}. (The estimate
$\|A_{k^*}^{-1}\|=O(1/\|k\|)$ mentioned in Remark
\ref{rmk:lemma-torus} is an easy consequence of
\eqref{eq:A_k-inverse}.)

The exact expression \eqref{eq:A_k-inverse-norm} can either be verified
by a direct calculation or proved as follows. Namely, using again the
fact that all complex structures $J_l$ anti-commute and are orthogonal
operators, it is easy to check that $(J\sum k_lJ_l)^{\top}=-J(\sum k_l
J_l)$. Then, from \eqref{eq:A_k} and \eqref{eq:A_k-inverse}, we infer
that $A_{k^*}$ and $A_{k^*}^{-1}$ are self-adjoint. Using again
\eqref{eq:A_k-inverse}, we have
$$
\left< A_{k^*}^{-1}Z, A_{k^*}^{-1}Z\right>
=\left<Z, A_{k^*}^{-1}A_{k^*}^{-1}Z\right>
=\frac{1}{4\pi^2\|k\|^2}\left<Z, A_{k^*}A_{k^*}^{-1}Z\right>
=\frac{1}{4\pi^2\|k\|^2}\|Z\|^2
$$
for any $Z\in F_{k*}$. This proves \eqref{eq:A_k-inverse-norm} and
completes the proof of the lemma.
\end{proof}

The rest of the argument, closely following \cite{CZ}, has become
quite standard by now and is included here only for the sake of
completeness. Denote by $\FF_N$ the subspace in $\FF$ formed by smooth
maps $f$ with $\hf_k=0$ whenever $\|k\|\geq N$. In other words, $\FF_N$
consists of Fourier polynomials of degree less than $N$, where the
degree is defined as $\|k\|$ in place of the more conventional
$|k|=|k_1|+\ldots+|k_r|$. Furthermore, let $\FF_N^\perp$ be the
$L^2$-orthogonal complement of $\FF_N$ in $\FF$, i.e., $\FF_N^\perp$
is the space of smooth maps $f$ with $\hf_k=0$ whenever $\|k\|< N$.  We
can view $\EE_N:=W\times \FF_N$ as a subbundle in $\EE$. 
It will also be useful to regard $\EE$ as a vector bundle over
$\EE_N$ with fiber $\FF_N^\perp$.  Denote by $\PP_N$ 
the (fiberwise) $L^2$-orthogonal projection of $\EE$ onto $\EE_N$ and
by $\PP_N^\perp$ the projection of $\EE=\EE_N\times\FF_N^\perp$ 
onto the second component $\FF_N^\perp$.

As is clear from Lemma \ref{lemma:torus}, the operator
$\partialslash|_{\FF_N^\perp}$ is invertible. Its inverse, which we
denote by $\partialslash_N^{-1}$, is $L^2$-bounded. Hence,
$\partialslash_N^{-1}$ extends by continuity to the $L^2$-completion
$\bar{\FF}_N^\perp$ of $\FF_N^\perp$. (The space $\bar{\FF}_N^\perp$
is formed by $L^2$-maps $f\colon M\to V$ with zero mean such that
$\hf_k=0$ for all $k$ with $\| k\|<N$.) Furthermore, again by
Lemma~\ref{lemma:torus}, we see that
\emph{$\|\partialslash_N^{-1}\|_{L^2}\leq 1/2\pi N$ and
  $\partialslash_N^{-1}$ sends functions of Sobolev class $H^s$ to
  functions of class $H^{s+1}$.} (The latter statement is, of
course, also a consequence of the fact, mentioned in Remark
\ref{rmk:ellipticity}, that $\partialslash$ is a first order elliptic
operator; see, e.g., \cite[Chap.\ III]{LM}.)

Our goal is to show that equation \eqref{eq:crit-pnts} has at least
the desired number of solutions. Let $f=g+h$ with $g\in \EE_N$ and
$h\in \FF_N^\perp$. Clearly, $f$ satisfies \eqref{eq:crit-pnts} if and
only if we have
\begin{equation}
\label{eq:crit-1}
\partialslash g=\PP_N \nabla H(g+h)
\end{equation}
and
\begin{equation}
\label{eq:crit-2}
\partialslash h =\PP_N^\perp \nabla H(g+h).
\end{equation}

Let us focus on the second of these equations with $g$ fixed and both
sides viewed as functions of $h$, cf.\ \cite{CZ}. Clearly,
\eqref{eq:crit-2} is equivalent to
\begin{equation}
\label{eq:crit-3}
h =\partialslash_N^{-1} \PP_N^\perp \nabla H(g+h).
\end{equation}
Note that the right hand side is now defined for all $h\in
\bar{\FF}_N^\perp$ without any smoothness requirement. We claim that
when $N$ is large enough, \emph{for any $g\in \EE_N$, equation
  \eqref{eq:crit-3} (and hence \eqref{eq:crit-2}) has a unique
  solution $h=h(g)$ and this solution is smooth}.

To show this, note first that, when $N$ is
sufficiently large, $h\mapsto \partialslash_N^{-1}
\PP_N^\perp \nabla H(g+h)$ is a contraction operator on
$\bar{\FF}_N^\perp$ with respect to the $L^2$-norm. Indeed,
$$
\|\partialslash_N^{-1} \PP_N^\perp \nabla H(g+h_1)
-\partialslash_N^{-1} \PP_N^\perp \nabla
H(g+h_0)\|_{L_2}
\leq \frac{1}{2\pi N}\| \nabla H (g+h_1)-\nabla H (g+h_0)\|_{L^2}
$$
and, in obvious notation,
\begin{eqnarray*}
\| \nabla H (g+h_1)-\nabla H (g+h_0)\|_{L^2} &=&
\Big\| \int_0^1 \frac{d}{d s} \nabla H (g+ s h_1 +(1-s)
h_0) \, ds \Big\|_{L^2}\\
&\leq& \|\nabla^2 H\|_{L^\infty}\|h_1-h_0\|_{L^2}.
\end{eqnarray*}
Hence,
$$
\|\partialslash_N^{-1} \PP_N^\perp \nabla H(g+h_1)-\partialslash_N^{-1} \PP_N^\perp \nabla
H(g+h_0)\|_{L_2}
\leq O(1/N) \|h_1-h_0\|_{L^2},
$$
which shows that we can indeed choose $N$ such that the map
$h\mapsto \partialslash_N^{-1} \PP_N^\perp \nabla H(g+h)$ is a
contraction. The fact that the fixed point $h=h(g)$ of this operator
is a smooth function is established by the standard
bootstrapping argument. Namely, we have $\partialslash h=\PP_N^\perp\nabla H (g+h)\in
L^2=H^0$, and therefore $h\in H^1$. Now, since $H$ and $g$ are smooth,
we also have $\PP_N^\perp\nabla H (g+h)\in H^1$, and hence $h\in H^2$, etc.

From a more geometrical perspective, $h(g)$ is the unique critical
point of the action functional $\A_H$ on the fiber over $g$ of the
vector bundle $\EE\to\EE_N$. Set $\Phi(g):=\A_H(g+h(g))$. In other
words, $\Phi$ is obtained from $\A_H$ by restricting the action functional
to the section $g\mapsto h(g)$ of this vector bundle, formed by the
fiber-wise critical points. Therefore, $g$ is a critical point of $\Phi$
if and only if $f=g+h(g)$ is a critical point of $\A_H$, i.e., a
solution of \eqref{eq:crit-pnts}, and every critical point of $\A_H$ is
captured in this way. It remains to show that the generating function
$\Phi$ on $\EE_N$ has the required number of critical points.

The key feature of this function is that it is asymptotically (i.e., at
infinity in the fibers of $\EE_N$) a non-degenerate quadratic
form. To be more precise, set
$$
\Phi_0(g)=\A(g)=\left<\partialslash g,g\right>_{L^2} \text{ and } R=\Phi-\Phi_0 .
$$
The unperturbed action $\Phi_0$ is a fiberwise non-degenerate quadratic
form. By definition, $\nabla \Phi_0(g)=\partialslash g$. (The quadratic form
$\Phi_0$ has zero signature, but this is not essential for what follows.)
Furthermore, the perturbation $R$ is small compared to $\Phi_0$, when $N$
is sufficiently large. Namely, for our purposes it is sufficient to show that
fiberwise
\begin{equation}
\label{eq:gen-function}
|R|+\|\nabla R\|<\| \nabla \Phi_0\| \text{ outside a compact set.}
\end{equation}
Here and throughout the rest of the proof, the metric on 
$\EE_N=W\times \FF_N$ is the product of the fiberwise $L^2$-metric and
the metric on $W$.

To establish \eqref{eq:gen-function}, note first that $H$ and $\nabla
H$ are bounded; for $H$ is a function on a compact manifold.
Therefore, the integral of $H$ makes a bounded contribution to $R$ and
$\nabla R$, while the right hand side of \eqref{eq:gen-function} grows
linearly as $g\to \infty$ in the fiber.  Thus, we can ignore $H$ in
\eqref{eq:gen-function} and only need to estimate the growth of the
difference
$$
R_0:=\A(g+h(g))-\A(g)=2\left<\partialslash g, h(g)\right>+
\left<\partialslash h(g), h(g)\right>,
$$
or to be more precise of $|R_0|$ together with $\| \nabla R_0\|$.
First observe that $| R_0(g) |$ is bounded by $O(1/N)(\|\nabla
\Phi_0(g)\|+1)$. (Here and below all the bounds are in the
$L^2$-norm.) This follows from the facts that the function $g\mapsto
h(g)$ is uniformly bounded by a constant $O(1/N)$, due to
\eqref{eq:crit-3}, and that the function $g\mapsto \partialslash h(g)$
is uniformly bounded, due to \eqref{eq:crit-2}. In a similar vein, it
is not hard to show that $\| \nabla R_0(g)\|$ is bounded from above by
$O(1) + O(1/N)\|\nabla \Phi_0 (g)\|$. (To this end, one also uses
the fact that the derivative of the function $g\mapsto h(g)$ is
uniformly bounded by a constant $O(1/N)$, as can be seen by
differentiating \eqref{eq:crit-3} with respect to $g$.)  Together,
these upper bounds prove~\eqref{eq:gen-function}.

A similar argument shows that a critical point $g$ of
$\Phi$ is non-degenerate when $f=g+h(g)$ is a non-degenerate critical
point of $\A_H$.

Finally, recall that whenever $\Phi=\Phi_0+R$ is a function on the total
space of a vector bundle over an arbitrary closed manifold $W$, such
that $\Phi_0$ is a fiberwise non-degenerate quadratic form and
\eqref{eq:gen-function} holds, the function $\Phi$ has at least
$\CL(W)+1$ critical points. Moreover, when $\Phi$ is Morse, the number of
critical points is bounded from below by $\SB(W)$. This is a
standard fact and we refer the reader to \cite{CZ} for the original
proof and to, e.g., \cite{We} for a different argument.  (Here we only
mention that the requirement \eqref{eq:gen-function} enables one to
modify $\Phi$ outside a sufficiently large compact set, without creating
new critical points, to turn it into a function identically equal to
$\Phi_0$ at infinity.)

Turning to the general case where $W$ is the quotient of $V$ by a
group $\Gamma$, we argue as follows. First recall that $\Gamma$
contains a finite-index subgroup $\Gamma'$ consisting of only parallel
transports, \cite[p. 110]{Wo}. Thus $W'=V/\Gamma'$ is a torus and the
projection $W'\to W$ is a covering map with the finite group
$\Pi=\Gamma/\Gamma'$ acting as the group of deck transformations. The
previous argument applies to the natural lift of the problem to $W'$
and the entire construction is $\Pi$-equivariant. As a result, we
obtain a vector bundle $\EE'_N\to W'$ equipped with a $\Pi$-action
covering the $\Pi$-action on $W'$ and a $\Pi$-invariant function $\Phi'$
on $\EE'_N$, which is asymptotically quadratic at infinity. The
critical points of $\A_H$ for the original problem correspond to the
$\Pi$-orbits of the critical points of $\Phi'$. Passing to the quotient
by $\Pi$, we arrive at a vector bundle over $W$ and a smooth function
$\Phi$ on its total space $\EE'_N/\Pi$. (The total space is smooth; for
the $\Pi$-action on $\EE'_N$ is free as an action covering a free
action on $W'$.) The function $\Phi$ is asymptotically quadratic and its
critical points are in one-to-one correspondence with the critical
points of $\A_H$ for the original problem. The theorem now follows as 
before from the lower bounds on the number of critical points of $\Phi$.

\subsection{The $\SU(2)$-case}
\label{sec:su(2)}
Let us now consider the setting where $M=\SU(2)$ and $r=3$ and $W$ is
the quotient, by a lattice, of a hyperk\"ahler vector space $V$ with
complex structures $J_1$, $J_2$ and $J_3$. (In particular, $W$ is a
torus.)  The case of a more general quotient $W=V/\Gamma$ can be
reduced to this one exactly as in Section \ref{sec:torus}; see the
previous paragraph. Furthermore, the case where $M$ is the quotient
$\SU(2)/G$ does not present any new difficulties and in fact follows
from the argument below. Throughout the rest of the proof, we will
treat $V$ as a real vector space or as a complex vector space with
complex structure $J=J_3$. Let us also fix a Hermitian
inner product on $V$, which, when necessary, we can also view as a
real inner product by discarding the imaginary part.

The space $L^2(\SU(2),V)$ is a unitary representation of $\SU(2)$,
which, by the Peter--Weyl theorem, decomposes into an $L^2$-sum of
irreducible representations $P_k$, $k=0,1,2,\ldots$, of $\SU(2)$ with
$P_k$ entering the sum with multiplicity $\dim_{\C} (P_k\otimes V)$;
see, e.g., \cite{Bo}.

The irreducible representation $P_k$ is the natural representation of
$\SU(2)$ on the space of homogeneous polynomials of degree $k$ in two
complex variables $z_1$ and $z_2$.  The $\SU(2)$-action on $P_k$ is
given by $x\cdot p=p\circ x^{-1}$ for $p\in P_k$ and $x\in\SU(2)$. Let
us turn $P_k$ into a unitary representation by fixing a Hermitian
inner product $\left<\cdot,\cdot\right>$ on $P_k$ which is invariant
under the group action. (Note that that such an inner product is
unique up to a factor; the normalization of the inner product is
immaterial for what follows.) Set $e_a^{(k)}=z_1^az_2^{k-a}$ for
$a=0,\ldots, k$. This is an orthogonal basis of $P_k$ with respect to
$\left<\cdot,\cdot\right>$.  The matrix coefficients
$e_{a,b}^{(k)}\colon \SU(2)\to \C$ for $a,b\in\{0,\ldots,k\}$ are
defined as
$$
e_{a,b}^{(k)}(x)=\big< x\cdot e_a^{(k)},\ e_b^{(k)}\big>.
$$
These are complex--valued functions on $\SU(2)$. With $i$ acting as $J=J_3$, 
we will view matrix coefficients as $\GL(V)$--valued functions.  

As in the torus case, the domain $\EE$ of the action functional $\A_H$
consists of smooth null-homotopic functions $f\colon \SU(2)\to W$.
Such a function $f$ lifts to an $L^2$-map $\tilde{f}\colon \SU(2)\to
V$.  Using the Peter--Weyl theorem, we can decompose $\tilde{f}$ as
$$
\tilde{f}(x) =\sum_{k\geq 0}\ 
\sum_{a,b=0}^k  e_{a,b}^{(k)}(x)\hat{f}_{a,b}^{(k)}.
$$
Here the sum converges in $L^2(\SU(2),V)$, the terms are mutually
$L^2$-orthogonal, and the Fourier coefficients $\hat{f}_{a,b}^{(k)}\in V$
are uniquely determined by $\tilde{f}$. (The same of course holds for
any $V$-valued $L^2$-function on $\SU(2)$.)

It is essential for what follows that
the coefficients of the non-constant matrix elements, i.e., the vectors
$\hat{f}_{a,b}^{(k)}\in V$ for $k\neq 0$, depend only on $f$ and are
independent of the lift.  As in the torus case, we can therefore
write, slightly abusing notation,
\begin{equation}
\label{eq:fourier-exp-su2}
f(x)=\sum_{k\geq 0}\ \sum_{a,b=0}^k  e_{a,b}^{(k)}(x) \hat{f}_{a,b}^{(k)}\ ,
\end{equation}
where $\hat{f}_{a,b}^{(k)}\in V$ when $k\neq 0$ and the mean value
$\hat{f}_{0,0}^{(0)}$ is an element of $W$.  Thus, the space $\EE$ can
be viewed as an infinite-dimensional vector bundle over $W$
with projection map $f\mapsto\hat{f}_{0,0}^{(0)}$. This vector bundle
is trivial and its fiber $\FF$ is canonically isomorphic to the space
of smooth maps $\SU(2)\to V$ with zero mean.

For a fixed $k>0$, we denote by $F_k$ the subspace of $L^2(\SU(2),V)$
which is spanned by all functions $e_{a,b}^{(k)}(x) w$ for $w\in V$ and
$a,b\in\{0,\ldots, k\}$. This is the subspace formed by the functions $f$
such that $\hat{f}_{a,b}^{(l)}=0$ for $l\neq k$.

We are now in a position to find an explicit representation of the
operator $\partialslash f$ in terms of the Fourier expansion of $f$.
The image of a function $f$ under $\partialslash$ is independent of
the mean value $\hat{f}_{0,0}^{(0)}$, since the constant term is
killed by the derivatives in $\partialslash$. Therefore, we can view
$\partialslash$ as a fiberwise linear map on $\EE=W\times \FF$, which
is independent of the point in the base $W$. Our goal is to
block-diagonalize $\partialslash$. In what follows, it is useful to
keep in mind that this operator is not complex linear.

In order to identify the invariant subspaces of $\partialslash$, we
utilize the decomposition of $\FF$ over irreducible representations
along with the quaternionic structure on $V$.  To be more precise,
recall that $V$ is not just a complex vector space, but also a
quaternionic vector space; for the complex structures $J_1,\, J_2,\,
J_3$ satisfy the quaternionic relations. Thus, we can decompose $V$ as
the sum of four real vector spaces intertwined by the operators $J_m$,
i.e., $V=V_0\oplus V_1\oplus V_2\oplus V_3=V_0^4$, where $V_m=J_mV_0$
for $m=1,2,3$. Let us denote by $I$ the identity map on $V$ or $V_0$.

\begin{Lemma}
\label{lemma:su2}
The operator $\partialslash$ preserves the subspaces $F_k$ and on each
of this subspaces block-diagonalizes as the sum of the following
operators:

\begin{enumerate}

\item the scalar operator $k\cdot I$ on the space $e_{a,0}^{(k)}\ V$,

\item the operator

\begin{equation}
\label{matrix0}
\left[
\begin{array}{cc}
(k-2b)I & (-1)^{a+b}2(k-b+1)I\\
(-1)^{a+b}2b I & (2b-k-2) I
\end{array}
\right]
\end{equation}
on the subspace $e_{a,b}^{(k)}\ V_0\ \oplus\
e_{k-a,k-b+1}^{(k)}V_2\cong V_0\oplus V_0$ for $b\in\{1,\ldots, k\}$,

\item the operator
\begin{equation}
\label{matrix1}
\left[
\begin{array}{cc}
(k-2b)I & (-1)^{a+b+1}2(k-b+1)I\\
(-1)^{a+b+1}2b I & (2b-k-2) I
\end{array}
\right]
\end{equation}
on the space $e_{a,b}^{(k)}\ V_1\ \oplus\ e_{k-a,k-b+1}^{(k)}V_3\cong
V_0\oplus V_0$ for $b\in\{1,\ldots, k\}$.

\end{enumerate}
Furthermore, on each of the subspaces $F_k$ with $k>0$, the operator
$\partialslash$ is invertible and its inverse has norm $1/k$.
\end{Lemma}

\begin{Remark}
\label{rmk:lemma-SU(2)} 
As in the torus case, this lemma is more precise than is really
 necessary for the proof. It would be sufficient to know that
$\partialslash|_{F_k}$ is invertible for $k>0$ and that its inverse has norm $O(1/k)$.
\end{Remark}

\begin{proof}
  First let us determine the matrix representation of the operator
  $\partialslash$ and show that the subspaces $F_k$ are invariant.
  Recall that the operator $\partialslash$ is given by
$$
\partialslash =J_1L_{v_1} + J_2L_{v_2} + J_3L_{v_3},
$$
where the right-invariant vector fields $v_l$ are defined by their
values at the identity as in \eqref{vectorfields}.  Computing the Lie
derivatives yields
\[
J_1\ L_{v_1}e_{a,b}^{(k)}=J_1J_3\left( b\ e_{a,b-1}^{(k)} + (k-b)\
  e_{a,b+1}^{(k)} \right),
\]
\[
J_2\ L_{v_2}e_{a,b}^{(k)}=J_2\left(- b\ e_{a,b-1}^{(k)} + (k-b)\
  e_{a,b+1}^{(k)}\right)
\]
and
\[
J_3\ L_{v_3}e_{a,b}^{(k)}=J_3\left(-J_3 (k-2b)\ e_{a,b}^{(k)}\right).
\]
Here we set $e_{a,-1}^{(k)}=0=e_{a,k+1}^{(k)}$. (In fact, the actual
definition of these functions is immaterial since they enter the
formulas with zero coefficients.)  Taking into account the
quaternionic relations between the complex structures, we obtain
\begin{equation}
\label{partialslash1}
\partialslash e_{a,b}^{(k)}=(k-2b)\ e_{a,b}^{(k)} - J_2\ 2b\ e_{a,b-1}^{(k)}.
\end{equation}

When $b=0$, this is the result of part (i) of
the lemma.

To deal with the case $b\neq 0$, recall first that $J_2$ is not a complex
linear operator on $V$: it does not commute with $J=J_3$. However, it
anti-commutes with $J$, i.e., $J_2J=-JJ_2$, and hence
$$
J_2\ e_{a,b}^{(k)} = \overline{e^{(k)}_{a,b}}\ J_2,
$$
since the matrix coefficients are $\C$-valued functions.
A direct calculation of the matrix coefficients or an
argument using the conjugate representation of $\SU(2)$ yields that
$$
\overline{e^{(k)}_{a,b}}\ J_2 = (-1)^{a+b} e_{k-a,k-b}^{(k)}\ J_2.
$$
Using this, we can rewrite \eqref{partialslash1} for $b\neq 0$ as
$$
\partialslash e_{a,b}^{(k)} =(k-2b)\ e_{a,b}^{(k)} + 
(-1)^{a+b} 2b\ e_{k-a,k-b+1}^{(k)}\ J_2
$$
With the identifications $V_0\oplus J_2V_0 = V_0\oplus
V_2\cong V_0\oplus V_0$ and $V_1 \oplus J_2V_1\cong V_1\oplus V_3\cong
V_0\oplus V_0$, this formula immediately implies the matrix
representations given in parts (ii) and (iii).

Let us now turn to the ``moreover" part of the lemma
and prove the bounds on the inverse of $\partialslash$. It is not hard
to check that the operator $\partialslash$ is invertible on $F_k$ for
$k>0$ by computing the eigenvalues of the matrices. In the case $b=0$
in part (i), it is clear that $k$ is the only eigenvalue. On the
subspaces considered in parts (ii) and (iii), one easily computes the
eigenvalues to be $k$ and $-k-2$. Thus, zero is not an eigenvalue for
$k>0$ and the operator is invertible and, moreover, the inverse of
$\partialslash|_{F_k}$, has eigenvalues $1/k$ and
$-1/(k+2)$. Furthermore, the eigenvectors are mutually orthogonal
since $\partialslash$ is self-adjoint. This shows that the norm of the
inverse of $\partialslash|_{F_k}$ is indeed $1/k$ as stated in the lemma.
\end{proof}

The remaining part of the proof of Theorem \ref{thm:main} goes through
almost word-for-word as in the torus case and in \cite{CZ}.  Recall
that we view the space of all null-homotopic functions as a trivial
vector bundle $\EE=W\times\FF$ and that the fiber $\FF$ is the direct
sum of the subspaces $F_k$ for $k>0$.  Closely following the reasoning
in the torus case (see Section \ref{sec:torus}), we denote the direct
sum of $F_k$ for $0<k< N$ by $\FF_N$ and its $L^2$-orthogonal
complement by $\FF_N^{\perp}$. Thus, $\FF_N$ consists of all functions
with $\hat{f}_{a,b}^{(k)}=0$ for $k\geq N$. For $f\in\FF_N^\perp$, we
have $\hat{f}_{a,b}^{(k)}=0$ for $k< N$.

Set $\EE_N=W\times\FF_N$ and
denote the fiberwise orthogonal projection of $\EE$ onto $\EE_N$ by
$\PP_N$. Let $\PP_N^\perp$ again denote the projection of 
$\EE=\EE_N\times\FF_N^\perp$ onto the second component.  
By Lemma \ref{lemma:su2}, the restriction
of the operator $\partialslash$ to $\FF_N^\perp$ is invertible and, on
$\FF_N^\perp$, the $L^2$-norm of the inverse
$\partialslash_N^{-1}:=(\partialslash|_{\FF_N^\perp})^{-1}$ is bounded
by $O(1/N)$. As a consequence, this operator extends by continuity to
the $L^2$-completion $\bar{\FF}_N^\perp$ of~$\FF_N^\perp$.

We need to show that equation \eqref{eq:crit-pnts} has at least the
desired number of solutions.  As in Section \ref{sec:torus}, we write $f=g+h$
with $g\in\EE_N$ and $h\in\FF_N^\perp$ and break the equation
\eqref{eq:crit-pnts} into equations \eqref{eq:crit-1} and
\eqref{eq:crit-2}.  For a fixed $g\in\EE_N$, equation
\eqref{eq:crit-2} gives rise to the fixed point problem
\begin{equation}
\label{eq:crit-4}
h =\partialslash_N^{-1} \PP_N^\perp \nabla H(g+h)
\end{equation}
for $h\in\FF_N^\perp$, where the right hand side is defined for all
$h\in \bar{\FF}_N^\perp$ without any smoothness requirement. We claim
that for any sufficiently large $N$ and any $g\in \EE_N$, equation
\eqref{eq:crit-4} has a unique solution $h=h(g)$ and that this
solution is smooth.

The existence of the solution $h(g)$ is established by the same
argument as in the torus case. Namely, the Hamiltonian $H$ is smooth
and compactly supported. Thus, $H$ and $\nabla H$ are uniformly
bounded by a constant. The norm of $\partialslash_N^{-1}$ is bounded 
by $O(1/N)$, due to Lemma \ref{lemma:su2}.  
For fixed $g$ and $H$, we can therefore
choose $N$ sufficiently large so that the operator
$h\mapsto\partialslash_N^{-1}\PP_N^\perp\nabla H(g+h)$ is a
contraction. This proves the existence and uniqueness of the fixed
point $h(g)$.
 
To show that $h(g)$ is smooth, we invoke elliptic regularity. Namely,
recall that, since $\partialslash$ is a first order elliptic operator
(see Remark \ref{rmk:ellipticity}), a solution $h$ of the equation
$\partialslash h= y$ is of Sobolev class $H^{s+1}$ whenever $h$ and
$y$ are $H^s$; see, e.g., \cite[Chap.\ III]{LM}. Applying this to $y
=\PP_N^\perp \nabla H(g+h)$ and using the standard bootstrapping
argument as in Section \ref{sec:torus}, we conclude that $h$ is
$C^\infty$-smooth.

From here on, the argument from the torus case applies without
any modifications. The calculations in Section \ref{sec:torus} are
independent of the specific setting of the torus case, relying only on
the definition of the action $\A_H$ by \eqref{eq:action}.  The
function $\Phi$ is asymptotically a non-degenerate quadratic form in
the fibers of the bundle $\EE\to\EE_N$. (Note that, in contrast to the
torus case, the quadratic form $\Phi$ does not have zero signature on
the subspaces $F_k$. However, this is not relevant for the proof of
the theorem.)  A critical point $g$ of $\Phi$ is non-degenerate if and
only if $f=g+h(g)$ is a non-degenerate critical point of $\A_H$.

Finally, recall that, as was already mentioned in Section
\ref{sec:torus}, a function $\Phi$ on the total space of a vector bundle
over a closed manifold $W$ has at least $\CL(W)+1$ critical points,
whenever $\Phi$ is asymptotically a non-degenerate quadratic form and $\Phi$
satisfies \eqref{eq:gen-function}. Moreover, when $\Phi$ is Morse, the
number of critical points is bounded from below by $\SB(W)$.  This
completes the proof of the theorem.


\begin{thebibliography}{BEHWZ}


\bibitem[Bo]{Bo}
N. Bourbaki,
\emph{Elements of Mathematics, Lie Groups and Lie Algebras},
Springer-Verlag, Berlin, Heidelberg, 2005.

\bibitem[CZ]{CZ}
C. Conley, E. Zehnder, 
The Birkhoff--Lewis fixed point theorem and a conjecture of V.I. Arnold,
\emph{Invent.\ Math.}, \textbf{73} (1983), 33--49.


\bibitem[FZ]{FZ} L. Florit, W. Ziller, Topological obstructions to
  fatness, Preprint 2010, arXiv:1001.0967.
 
\bibitem[GHJ]{GHJ} M. Gross, D. Huybrechts, D. Joyce, \emph{Calabi--Yau
    Manifolds and Related Geometries}, Lectures at a Summer School in
  Nordfjordeid, Norway, June 2001, Springer-Verlag, Berlin, Heidelberg,
  2003.
 
\bibitem[HNS]{HNS} S. Hohloch, G. Noetzel, D. Salamon, Hypercontact
  structures and Floer theory, \emph{Geom.\ Topol.}, \textbf{13}
  (2009), 2543--2617.

\bibitem[Hu]{Hu}
D. Husemoller, \emph{Fibre Bundles}, Springer-Verlag, New York, 1994.

\bibitem[Jo]{Jo}
D. Joyce, Manifolds with many complex structures, \emph{Quart.\
  J. Math.\ Oxford Ser.\ (2)}, \textbf{46} (1995), 169--184.

\bibitem[LM]{LM} H.B. Lawson, M.-L. Michelsohn, \emph{Spin Geometry},
  Princeton Mathematical Series, 38. Princeton University Press,
  Princeton, NJ, 1989.

\bibitem[Mo]{Mo} R. Montgomery, \emph{A Tour of Subriemannian
    Geometries, their Geodesics and Applications}, Mathematical
  Surveys and Monographs, vol.\ 91. American Mathematical Society,
  Providence, RI, 2002.

\bibitem[MS]{MS}
A. Moroianu, U. Semmelmann,
Clifford structures on Riemannian manifolds, Preprint 2009,
arXiv:0912.4207.

\bibitem[We1]{We:fat} A. Weinstein, Fat bundles and symplectic
  manifolds, \emph{Adv.\ in Math.}, \textbf{37} (1980), 239--250.

\bibitem[We2]{We}
A. Weinstein,  
$C^{0}$ perturbation theorems for symplectic fixed
points and Lagrangian intersections, in \emph{South Rhone seminar on
geometry, III (Lyon, 1983)}, pp.\ 140--144, \emph{Travaux en Cours}, 
Hermann, Paris, 1984.

\bibitem[Wo]{Wo}
J. Wolf,
\emph{Spaces of Constant Curvature}, Fifth edition,
Publish or Perish, Inc., Houston, TX, 1984

\end{thebibliography}
\end{document}